\documentclass{article}
\usepackage{amsfonts}

\usepackage{amsmath}

\newtheorem{Theorem}{Theorem}[section]
\newtheorem{Definition}[Theorem]{Definition}
\newtheorem{Proposition}[Theorem]{Proposition}
\newtheorem{Lemma}[Theorem]{Lemma}
\newtheorem{Corollary}[Theorem]{Corollary}
\newtheorem{Remark}[Theorem]{Remark}
\newtheorem{Example}[Theorem]{Example}

\def\qed{\hfill\hbox{\hskip 6pt\vrule width6pt height7pt
depth1pt  \hskip1pt}\bigskip}

\def\E{{{\rm I} \kern -.15em {\rm E}    }}

\newcommand{\reals}{{{\rm I} \kern -.15em {\rm R} }}
\newcommand{\complex}{{{\rm I} \kern -.52em {\rm C} }}
\newcommand{\nat}{{{\rm I} \kern -.15em {\rm N} }}

\def\endproof{\qquad \vbox to 5.8pt{\offinterlineskip\hrule
        \hbox to 5.8pt{\vrule height 5.1pt\hss\vrule height 5.1pt}\hrule}}

\begin{document}

\date{\today}

\vspace{10mm}

\begin{center}
\textbf{\Large Weak Dirichlet processes with a stochastic control
perspective.}

\vspace{ 5 mm}

{\large Fausto GOZZI}

\vspace{ 3 mm}

{\large Dipartimento di Scienze Economiche e Aziendali}\\[0pt]
{\large Facolta' di Economia }\\[0pt]
{\large LUISS - Guido Carli }\\[0pt]
{\large Viale Pola 12, I-00198  Roma, Italy }


\bigskip
{\large Francesco RUSSO}
\vspace{ 3 mm}

{\large Universit\'{e} Paris 13}\\[0pt]
{\large Institut Galil\'{e}e, Math\'{e}matiques}\\[0pt]
{\large 99, av. JB Cl\'{e}ment, F-99430 Villetaneuse, France}
\end{center}

\vspace{ 8 mm}

\textbf{A.M.S. Subject Classification}: 60G05, 60G44, 60G48, 60H05, 60H10,
60J60, 35K15, 35K55, 35J15, 93E20.


\vspace{ 3 mm}

\textbf{Key words}: 
Stochastic calculus via regularization, weak Dirichlet processes, stochastic
optimal control, Cauchy problem for parabolic partial differential equations.

\vspace{3mm}

\textbf{Abstract.} The motivation of this paper is to prove verification
theorems for stochastic optimal control of finite dimensional diffusion
processes without control in the diffusion term, in the case that the value
function is assumed to be continuous in time and once differentiable in
the space variable ($C^{0,1}$) instead of once differentiable in time
and twice in space ($C^{1,2}$), like in the classical results. For this
purpose, the replacement tool of the It\^{o} formula will be the
Fukushima-Dirichlet decomposition for weak Dirichlet processes. Given a
fixed filtration, a weak Dirichlet process is the sum of a local martingale $%
M$ plus an adapted process $A$ which is orthogonal, in the sense of
covariation, to any continuous local martingale.
The mentioned decomposition states
that a $C^{0,1}$ function of a weak Dirichlet process with finite quadratic
variation is again a weak Dirichlet process. That result is established in
this paper and it is applied to the strong solution of a Cauchy problem with
final condition.

Applications to the proof of verification theorems will be addressed in a
companion paper.

\section{Introduction\label{INTRO}}

In this paper we prepare a framework of stochastic calculus via
regularization in order to apply it to the proof of verification theorems in
stochastic optimal control in finite dimension. The application part will be
implemented in the companion paper \cite{GR2}.

This paper has an interest in itself and its most significant result is a
generalized time-dependent Fukushima-Dirichlet decomposition which is proved
in section \ref{ITO}. This will be the major tool for applications.

The proof of verification theorems for stochastic control problems under
classical conditions is an application of It\^{o} formula. In fact, under
good assumptions, the value function $V:\left[ 0,T\right] \times \mathbb{R}%
^{n}\rightarrow \mathbb{R}$ associated with a stochastic control problem is
of class $C^{1}$ in time and $C^{2}$ in space ($C^{1,2}$ in symbols). This
allows to apply it to the solution of a corresponding state equation $%
(S_{t}) $ and differentiate $V(t,S_{t})$ through the classical It\^{o}
formula, see e.g. \cite[pp. 140, 163, 172]{FS}. The substitution tool of
that formula will be a time-dependent Fukushima - Dirichlet decomposition
which will hold for functions $u:\left[ 0,T\right] \times \mathbb{R}%
^{n}\rightarrow \mathbb{R}$ that are $C^{0,1}$ in symbols; so, our
verification theorem will have the advantage of requiring less regularity on
the value function $V$ than the classical ones.

It is also possible to prove a verification theorem in the case when $V$ is
only continuous (see e.g. \cite{LYZ}, \cite[Section 5.2]{Zhoubook}, \cite
{GSZstochver}) in the framework of viscosity solutions: however such result
applied to our cases is weaker than ours, as it requires more assumptions on
the candidate optimal strategy; see  on this the last section of the
companion paper \cite{GR2} where also a comparison with other nonsmooth
verification theorems is performed.

We come back to the Fukushima-Dirichlet decomposition as replacement of
It\^o formula. Roughly speaking, given a function $u$ of class $C^{1,2}$,
classical It\^o formula gives a decomposition of $u(t, S_t) $ in a
martingale part, say $M$ (which is thrown away taking expectation in the
case of deterministic data and expected cost) plus an absolutely continuous
process, say $A$. Then, in case of deterministic data and expected cost, one
uses the fact that $u$ is a classical solution of a partial differential
equation (PDE), which is in fact the Hamilton-Jacobi-Bellman (HJB) equation,
to represent $A$ in term of the Hamiltonian function. If one wants to repeat
the above arguments when $u$ is not $C^{1,2}$, a natural way is to try to
extend the decomposition $u\left( \cdot ,S\right) = M + A$ and the
representation of $A$ via the HJB PDE. This is what we perform in this paper
in the case when $u\in C^{0,1}$ using a point of view that can be also
applied to problems with stochastic data and pathwise cost, so the HJB
becomes a stochastic PDE, see e.g. \cite{DPDMT,LS1,LS2,Peng}).

We propose in fact an extension of the classical Fukushima-Dirichlet
decomposition. That decomposition is inspired by the theory of Dirichlet
forms. A classical monography concerning this theory is \cite{Fu} where one
can find classical references on the subject. In Theorem 5.2.2 of \cite{Fu},
given a ``good''\ symmetric Markov process $\left( X_{t}^{x}\right) _{t\geq
0}$ and a function belonging to some suitable space (Dirichlet space), it is
possible to write
\begin{equation}
u\left( X_{t}\right) =u\left( x\right) +M_{t}^{u}+A_{t}^{u},
\label{eq:introFR}
\end{equation}
where $M$ is a local martingale and $\left( A_{t}^{u}\right) _{t\geq
0}$ is a zero quadratic variation process, for quasi-everywhere $x$,
i.e. for\ $x$\ belonging to a zero capacity set. For instance, if
$X=W$ is a classical Brownian motion in $\mathbb{R}^{n}$ then the
Dirichlet space is $H^{1}\left( \mathbb{R}^{n}\right) $ and
$M^{u}_t=\int_{0}^{t}\nabla u(W_{u})dW_{u}$. We call
(\ref{eq:introFR}) a Fukushima-Dirichlet decomposition. Our point of
view is of pathwise nature: as in \cite{fo, fod}, a process $Y$, as
$u(X)$, which is the sum of a local martingale and a zero quadratic
variation process, even without any link to Dirichlet forms, is
called a Dirichlet process. However, in \cite{fo}, the notion of
quadratic variation, even if in the same spirit as ours, was defined
through discretizations, while we define it through regularizations.
Therefore that notion of Dirichlet process is similar but not
identical to ours.

The papers \cite{RV3} and \cite{fp} reinterpret $A$ in a ''pathwise'' way as the
covariation process $\left[ \nabla u\left( B\right) ,B\right] $ transforming
(\ref{eq:introFR}) in a true It\^{o}'s formula; the first work considers $%
u\in C^{1}\left( \mathbb{R}^{n}\right) $ and it extends the framework to
reversible continuous semimartingales; the second work is connected with
Brownian motion and $u\in H^{1}\left( \mathbb{R}^{n}\right) $. The
literature on It\^{o}'s formula for non-smooth functions of semimartingales
or diffusion processes has known a lot of development in the recent years,
see for instance \cite{MN} for non-degenerate Brownian martingales, \cite
{frw2} for non-degenerate 1-dimensional diffusions with bounded measurable
drift or \cite{ERV} in the jump case.

In our applications, the fact of identifying the remainder process $A^{u}$
as a covariation is not so important since the goal is to give the
representation of it via the data of the HJB PDE. So we come back to the
spirit of the Fukushima-Dirichlet decomposition. Besides our ``pathwise''
approach to Dirichlet processes, the true novelty of this approach is the
time-inhomogeneous version of the decomposition; this is in particular
motivated by non-autonomous problems in control theory.

This is based on the theory, under construction, of weak Dirichlet processes
with respect to some fixed filtration $(\mathcal{F}_{t})$. A weak Dirichlet
process is the sum of a local martingale $M$ and a process $A$ which is
adapted and $[A,N]=0$ for any continuous
$(\mathcal{F}_{t})$-local martingale $N$. We will
be able in particular to decompose $u(t,D_{t})$ when $u\in C^{0,1}$ and $D$
a weak Dirichlet process with finite quadratic variation process, so in
particular if $D$ is a semimartingale (even diffusion process). This will be
our time-dependent Fukushima-Dirichlet decomposition: it will be the object
of Proposition \ref{pr:FDdec2} and Corollary \ref{Co:FDdec2}. In particular
that result holds for semimartingales (and so for diffusion processes).
The notion of weak Dirichlet process appears also in \cite{ERbis}. Our
Fukushima-Dirichlet decomposition could be linked to the theory of
``time-dependent Dirichlet forms'' developed for instance by \cite
{Os,Stan,Trut} but we have not  investigated that direction.

The paper is organized as follows. In section \ref{PRELIMINARIES} we
introduce some notations on real analysis and we establish preliminary
notions on calculus via regularization with some remarks on classical
Dirichlet processes. Section 3 will be devoted to some basic facts about
weak Dirichlet process and to the above mentioned Fukushima-Dirichlet
decomposition of process $(u(t,D_{t}))$, with some sufficient condition to
guarantee that the resulting process is a true Dirichlet process. Section
\ref{REPRESENTATION} will be concerned with application to the case where $u$
is a strong $C^{0,1}$ solution of a Cauchy parabolic problem with initial
condition; $C^{1}$ solutions of an elliptic problem are also
represented probabilistically.

\section{Preliminaries\label{PRELIMINARIES}}

\subsection{Notation\label{NOTATIONS}}

Throughout this paper we will denote by $\left( \Omega ,\mathcal{F},\mathbb{P%
}\right) $ a given stochastic basis (where $\mathcal{F}$ stands for a given
filtration $\left( \mathcal{F}_{s}\right) _{s\geq 0}$ satisfying the usual
conditions). Given a finite dimensional real Hilbert space $E$, $W$ will
denote a cylindrical Brownian motion with values in $E$ and adapted to $%
\left( \mathcal{F}_{s}\right) _{s\geq 0}$. Given $0\leq t\leq T\leq +\infty $
and setting $\mathcal{T}_{t}=\left[ t,T\right] \cap \mathbb{R}$ the symbol $%
\mathcal{C}_{\mathcal{F}}\left( \mathcal{T}_{t}\times \Omega ;E\right) $,
will denote the space of all continuous processes adapted to the filtration $%
\mathcal{F}$ with values in $E$. This is a Fr\'{e}chet space if endowed with
the topology of the uniform convergence in probability ($u.c.p.$ from now
on). To be more precise this means that, given a sequence $\left(
X^{n}\right) \subseteq \mathcal{C}_{\mathcal{F}}\left( \mathcal{T}_{t}\times
\Omega ;E\right) $ and $X\in \mathcal{C}_{\mathcal{F}}\left( \mathcal{T}%
_{t}\times \Omega ;E\right) $ we have
\begin{equation*}
X^{n}\rightarrow X
\end{equation*}
if and only if for every $\varepsilon >0$, $t_{1}\in \mathcal{T}_{t}$%
\begin{equation*}
\lim_{n\rightarrow +\infty }\sup_{s\in \left[ t,t_{1}\right] }\mathbb{P}%
\left( \left| X_{s}^{n}-X_{s}\right| _{E}>\varepsilon \right) =0.
\end{equation*}

Given a random time $\tau \geq t$ and a process $(X_{s})_{s\in \mathcal{T}%
_{t}},$ we denote by $X^{\tau }$ the stopped process defined by $X_{s}^{\tau
}=X_{s\wedge \tau }.$ The space of all processes in $\left[ t,T\right] $,
adapted to $\mathcal{F}$ and square integrable with values in $E$ is denoted
by $L_{\mathcal{F}}^{2}\left( t,T;E\right) $. $S^{n}$ will denote the space
of all symmetric matrices of dimension $n$. If $Z$ is a vector or a matrix,
then $Z^{\ast }$ is its transposition.

Let $k\in \mathbb{N}$. As usual $C^{k}\left( \mathbb{R}^{n}\right) $ is the
space of all functions $:\mathbb{R}^{n}\rightarrow \mathbb{R}$ that are
continuous together with their derivatives up to the order $k$. This is a
Fr\'{e}chet space equipped with the seminorms
\begin{equation*}
\sup_{x\in K}\left| u\left( x\right) \right| _{\mathbb{R}}+\sup_{x\in
K}\left| \partial _{x}u\left( x\right) \right| _{\mathbb{R}^{n}}+\sup_{x\in
K}\left| \partial _{xx}u\left( x\right) \right| _{\mathbb{R}^{n\times n}}
\end{equation*}
for every compact set $K\subset \subset \mathbb{R}^{n}$. This space will be
denoted simply by $C^{k}$ when no confusion may arise. The symbol $%
C_{b}^{k}\left( \mathbb{R}^{n}\right) $ will denote the Banach space of all
continuous and bounded functions from $\mathbb{R}^{n}$ to $\mathbb{R}$. This
space is endowed with the usual $\sup $ norm. Passing to parabolic spaces we
denote by $C^{0}\left( \mathcal{T}_{t}\times \mathbb{R}^{n}\right) $ the
space of all functions
\begin{equation*}
u:\mathcal{T}_{t}\times \mathbb{R}^{n}\rightarrow \mathbb{R},\quad \quad
\left( s,x\right) \mapsto u\left( s,x\right)
\end{equation*}
that are continuous. This space is a Fr\'{e}chet space equipped with the
seminorms
\begin{equation*}
\sup_{\left( s,x\right) \in \left[ t,t_{1}\right] \times K}\left| u\left(
s,x\right) \right| _{\mathbb{R}}
\end{equation*}
for every $t_{1}>0$ and every compact set $K\subset \subset \mathbb{R}^{n}$%
). Moreover we will denote by $C^{1,2}\left( \mathcal{T}_{t}\times \mathbb{R}%
^{n}\right) $ (respectively $C^{0,1}\left( \mathcal{T}_{t}\times \mathbb{R}%
^{n}\right) $), the space of all functions
\begin{equation*}
u:\mathcal{T}_{t}\times \mathbb{R}^{n}\rightarrow \mathbb{R},\quad \quad
\left( s,x\right) \mapsto u\left( s,x\right)
\end{equation*}
that are continuous together with their derivatives $\partial _{t}u$, $%
\partial _{x}u$, $\partial _{xx}u$ (respectively $\partial _{x}u$). This
space is a Fr\'{e}chet space equipped with the seminorms
\begin{eqnarray*}
&&\sup_{\left( s,x\right) \in \left[ t,t_{1}\right] \times K}\left| u\left(
s,x\right) \right| _{\mathbb{R}}+\sup_{\left( s,x\right) \in \left[ t,t_{1}%
\right] \times K}\left| \partial _{s}u\left( s,x\right) \right| _{\mathbb{R}%
^{n}} \\
&&+\sup_{\left( s,x\right) \in \left[ t,t_{1}\right] \times K}\left|
\partial _{x}u\left( s,x\right) \right| _{\mathbb{R}^{n}}+\sup_{\left(
s,x\right) \in \left[ t,t_{1}\right] \times K}\left| \partial _{xx}u\left(
s,x\right) \right| _{\mathbb{R}^{n\times n}}
\end{eqnarray*}
(respectively
\begin{equation*}
\sup_{\left( s,x\right) \in \left[ t,t_{1}\right] \times K}\left| u\left(
s,x\right) \right| _{\mathbb{R}}+\sup_{\left( s,x\right) \in \left[ t,t_{1}%
\right] \times K}\left| \partial _{x}u\left( s,x\right) \right| _{\mathbb{R}%
^{n}}\text{)}
\end{equation*}
for every $t_{1}>0$ and every compact set $K\subset \subset \mathbb{R}^{n}$.
This space will be denoted simply by $C^{1,2}$ (respectively $C^{0,1}$) when
no confusion may arise.

Similarly, for $\alpha ,\beta \in \left[ 0,1\right] $ one defines $C^{\alpha
,1+\beta }\left( \mathcal{T}_{t}\times \mathbb{R}^{n}\right) $ (or simply $%
C^{\alpha ,1+\beta }$) as the subspace of $C^{0,1}\left( \mathcal{T}%
_{t}\times \mathbb{R}^{n}\right) $ of functions $u:\mathcal{T}_{t}\times
\mathbb{R}^{n}\mapsto \mathbb{R}$ such that are $u\left( \cdot ,x\right) $
is $\alpha -$H\"{o}lder continuous and $\partial _{x}u\left( s,\cdot \right)
$ is $\beta -$H\"{o}lder continuous (with the agreement that $0$-H\"older
continuity means just continuity).

Similarly to $C_{b}^{k}\left( \mathbb{R}^{n}\right) $ we define the Banach
spaces $C_{b}^{0}\left( \mathcal{T}_{t}\times \mathbb{R}^{n}\right) $ $%
C_{b}^{1,2}\left( \mathcal{T}_{t}\times \mathbb{R}^{n}\right) $, $%
C_{b}^{\alpha ,1+\beta }\left( \mathcal{T}_{t}\times \mathbb{R}^{n}\right) $%
, $C_{b}^{0,1}\left( \mathcal{T}_{t}\times \mathbb{R}^{n}\right) $.

\subsection{The calculus via regularization \label{RUSSOCALCULUS}}

We will follow here a framework of calculus via regularizations started in
\cite{RV1}. At the moment many authors have contributed to it and we suggest
the reader to consult the recent survey paper \cite{RVSem} on it.

For simplicity, all the considered processes, excepted if we mention the
contrary, will be continuous processes. We first recall some one dimensional
consideration. For two processes $\left( X_{s}\right) _{s\ge 0}$, $\left(
Y_{s}\right) _{s\ge 0}$, we define the forward integral and the covariation
as follows

\begin{equation}
\int_{0}^{s}X_{r}d^{-}Y_{r}=\lim_{\varepsilon \rightarrow 0}\int_{0}^{s}X_{r}%
\frac{Y_{r+\varepsilon }-Y_{r}}{\varepsilon }dr,  \label{eq:int-}
\end{equation}
\begin{equation}
\left[ X,Y\right] _{s}=\lim_{\varepsilon \rightarrow 0}\frac{1}{\varepsilon }%
\int_{0}^{s}\left( X_{r+\varepsilon }-X_{r}\right) \left( Y_{r+\varepsilon
}-Y_{r}\right) dr,  \label{eq:covdef}
\end{equation}
if those quantities exist in the sense of $u.c.p$ with respect to $s$. This
ensures that the forward integral defined in (\ref{eq:int-}) and the
covariation process defined in (\ref{eq:covdef}) are continuous processes.
It can be seen that the covariation is a bilinear and symmetric operator.

We fix now, as above, $0\leq t\leq T\leq +\infty $ and set $\mathcal{T}_{t}=%
\left[ t,T\right] \cap \mathbb{R}$. A process $\left( X_{s}\right)
_{s\in \mathcal{T}_{t}}$ can always be extended (if $T<+\infty $) to
a process indexed by $\mathbb{R}_{+}$ by continuity. The
corresponding extension will always denoted by the same symbol.
Given two processes $\left( X_{s}\right) _{s\in \mathcal{T}_{t}}$,
$\left( Y_{s}\right) _{s\in \mathcal{T}_{t}}$, we define the
corresponding stochastic integrals and covariations by the integrals
and covariation of the corresponding extensions. We define also
integrals from $t$ to $s$ as follows.
\begin{equation*}
\int_{t}^{s}X_{r}d^{-}Y_{r}=\int_{0}^{s}X_{r}d^{-}Y_{r}-%
\int_{0}^{t}X_{r}d^{-}Y_{r}.
\end{equation*}
If $\tau \geq t$ is a random (non necessarily) stopping time, the following
equality holds:
\begin{equation}
\lbrack X,Y]^{\tau }=[X^{\tau },Y^{\tau }].  \label{eq:bracket}
\end{equation}
If $\left( X^{1},...,X^{n}\right) $ is a vector of continuous processes we
say that it has all its mutual covariations (brackets) if $\left[ X^{i},X^{j}%
\right] $ exist for any $1\leq i,j\leq n$. If $X^{1},...,X^{n}$ have all
their mutual covariations then by polarisation (i.e. writing a bilinear form
as a sum/difference of quadratic forms) we know that $\left[ X^{i},X^{j}%
\right] $ ($1\leq i,j\leq n$) are locally bounded variation processes.

If $\left[ X,X\right] $ exists, then $X$ is said to be a finite quadratic
variation process; $\left[ X,X\right] $ is called the quadratic variation of
$X$. If $\left[ X,X\right] =0$ then $X$ is said to be a zero quadratic
variation process. A bounded variation process is a zero quadratic variation
process. If $S^{1}$, $S^{2}$ are $\left( \mathcal{F}_{s}\right) $%
-semimartingales then $\left[ S^{1},S^{2}\right] $ coincides with the
classical bracket $\left\langle S^{1},S^{2}\right\rangle $. If $H$ is a $%
\left( \mathcal{F}_{s}\right) $-progressively measurable process then $%
\int_{t}^{s}H_{r}d^{-}S_{r}$ is the classical It\^{o} integral $%
\int_{t}^{s}H_{r}dS_{r}$.

\begin{Remark}
\label{rm:bracketnullo}\textrm{Let }$X$\textrm{\ (respectively }$A$\textrm{)
be a finite (respectively zero) quadratic variation process. Then }$\left(
X,A\right) $\textrm{\ has all its mutual covariations and }$\left[ X,A\right]
=0$\textrm{.}\hfill
\hbox{\hskip 6pt\vrule width6pt height7pt
depth1pt  \hskip1pt}\bigskip
\end{Remark}

We recall now an easy extension of stability results, see \cite[th. 2.9]{FR}
that will be used in subsection \ref{DECHOLDER}.

\begin{Proposition}
\label{pr:bracket} Let $V=\left( V^{1},...V^{m}\right) $ (respectively $%
X=\left( X^{1},...X^{n}\right) $) be a vector of continuous processes on $%
\mathbb{R}_{+}$ with bounded variation processes (respectively having all
its mutual covariations). Let $f,g\in C_{\mathrm{loc}}^{\frac{1}{2}+\gamma
,1}\left( \mathbb{R}^{m}\times \mathbb{R}^{n}\right) $ ($\gamma >0$). Then $%
\forall s\geq 0$%
\begin{equation*}
\left[ f\left( V,X\right) ,g\left( V,X\right) \right] _{s}=\sum_{i,j=1}^{n}%
\int_{0}^{s}\partial _{x_{i}}f\left( V,X\right) \partial _{x_{j}}g\left(
V,X\right) d\left[ X^{i},X^{j}\right] _{r}.
\end{equation*}
\end{Proposition}

\textbf{Proof.}
We give a sketch of the proof as the arguments are similar to the ones used
in \cite[th. 2.9]{FR}. By localisation $C_{\mathrm{loc}}^{\frac{1}{2}+\gamma
,1}$ can be replaced by $C^{\frac{1}{2}+\gamma ,1}$. The case where $f$ and $%
g$ do not depend on $V$ was treated for instance in \cite{RV2,RV4}. Since
the covariation is a bilinear operation, using polarisation techniques we
can take $g=f$. For simplicity we set here $m=n=1$. For given $\varepsilon
>0 $ we write, when $r\in \mathcal{T}_{t}$,
\begin{equation*}
f\left( V_{r+\varepsilon },X_{r+\varepsilon }\right) -f\left(
V_{r},X_{r}\right) = J_{1}\left( r,\varepsilon \right) +J_{2}\left(
r,\varepsilon \right) ,
\end{equation*}
where
\begin{eqnarray*}
J_{1}\left( r,\varepsilon \right) &=&f\left( V_{r+\varepsilon
},X_{r+\varepsilon }\right) -f\left( V_{r},X_{r+\varepsilon }\right) , \\
J_{2}\left( r,\varepsilon \right) &=&f\left( V_{r},X_{r+\varepsilon }\right)
-f\left( V_{r},X_{r}\right) .
\end{eqnarray*}
Therefore, for $s\geq 0$ 
\begin{equation*}
\frac{1}{\varepsilon }\int_{0}^{s}\left[ f\left( V_{r+\varepsilon
},X_{r+\varepsilon }\right) -f\left( V_{r},X_{r}\right) \right] ^{2}dr\leq
\frac{2}{\varepsilon }\int_{0}^{s}J_{1}^{2}\left( r,\varepsilon \right) dr+%
\frac{2}{\varepsilon }\int_{0}^{s}J_{2}^{2}\left( r,\varepsilon \right) dr.
\end{equation*}
Now
\begin{equation*}
\frac{2}{\varepsilon }\int_{0}^{s}J_{1}^{2}\left( r,\varepsilon \right)
dr\leq c^{2}\left( s,f,V\right) \frac{2}{\varepsilon }\int_{0}^{s}\left(
\left| V_{r+\varepsilon }-V_{r}\right| ^{\frac{1}{2}+\gamma }\right) ^{2}dr,
\end{equation*}
where $c\left( s,f,V\right) $ is a (random) H\"{o}lder constant of $f$ such
that
\begin{equation*}
\left| f\left( v_{1},x\right) -f\left( v_{2},x\right) \right| \leq c\left(
s,f,V\right) \left| v_{1}-v_{2}\right| ^{\frac{1}{2}+\gamma },
\end{equation*}
\begin{equation*}
\forall v_{1},v_{2}\in \left[ \inf_{r\in \left[ 0,s\right] }V_{r},\sup_{r\in %
\left[ 0,s\right] }V_{r}\right] ,\quad \forall x\in \left[ \inf_{r\in \left[
0,s\right] }X_{r},\sup_{r\in \left[ 0,s\right] }X_{r}\right] .
\end{equation*}
Then we get
\begin{equation*}
\frac{2}{\varepsilon }\int_{0}^{s}J_{1}^{2}\left( r,\varepsilon \right)
dr\leq c^{2}\left( s,f,V\right) \frac{2}{\varepsilon }\int_{0}^{s}\left|
V_{r+\varepsilon }-V_{r}\right| ^{1+2\gamma }dr.
\end{equation*}
Since $V$ is a bounded variation process, this term converges to zero in
probability.

On the other hand,
\begin{equation*}
J_{2}\left( r,\varepsilon \right) =\partial _{x}f\left( V_{r},X_{r}\right)
\left( X_{r+\varepsilon }-X_{r}\right) +J_{3}\left( r,\varepsilon \right) ,
\end{equation*}
where $J_{3}\left( r,\varepsilon \right) $ converges $u.c.p.$ to zero, as in
\cite[th. 2.9]{FR}. Therefore, similarly as in \cite{RV2} we have
\begin{equation*}
\frac{2}{\varepsilon }\int_{0}^{t}J_{2}^{2}\left( s,\varepsilon \right)
ds\rightarrow \sum_{i,j=1}^{n}\int_{0}^{t}\left( \partial _{x}f\left(
V_{s},X_{s}\right) \right) ^{2}d\left[ X,X\right] _{s}
\end{equation*}
and so the claim.\hfill
\hbox{\hskip 6pt\vrule width6pt height7pt
depth1pt  \hskip1pt}\bigskip

For our purposes we need to express integrals and covariation in a
multidimensional setting, in the spirit of \cite{ERbis}.

If $X=\left( X^{1},...,X^{n}\right) ^{\ast }$ is a vector of continuous
processes in $\mathbb{R}_{+}$, $Y$ is a $m\times n$ matrix of continuous
processes in $\mathbb{R}_{+}$, $\left( Y^{i,j}\right) _{1\leq i\leq m,1\leq
j\leq n}$ then the symbol $\int_{0}^{s}Yd^{-}X$ denotes, whenever it exists,
the $u.c.p.$ limit of the integral $\int_{0}^{s}Y_{r}\frac{X_{r+\varepsilon
}-X_{r}}{\varepsilon }dr$ where the product is intended in the matrix sense.
Similarly, il $A$ is a $n\times d$ matrix $\left( A^{j,k}\right) _{1\leq
j\leq n,1\leq k\leq d}$ then $\left[ Y,A\right] _{s}$ is the $m\times d$
real matrix constituted by the following $u.c.p.$ limit (if it exists)
\begin{equation*}
\frac{1}{\varepsilon }\int_{0}^{s}\left( Y_{r+\varepsilon }-Y_{r}\right)
\left( A_{r+\varepsilon }-A_{r}\right) dr.
\end{equation*}
Clearly the matrix operation cannot be commutative in general. Let now $%
A,X,Y,C$ be real matrix valued processes which are successively compatible
for the matrix product. We define
\begin{equation*}
\int_{0}^{s}A_{r}d\left[ X,Y\right] _{r}C_{r}=\lim_{\varepsilon \rightarrow
0}\frac{1}{\varepsilon }\int_{0}^{s}A_{r}\left( X_{r+\varepsilon
}-X_{r}\right) \left( Y_{r+\varepsilon }-Y_{r}\right) C_{r}dr,
\end{equation*}
where the limit is intended in the $u.c.p.$ sense. The previous stability
transformations (Proposition \ref{pr:bracket} above) can be extended to the
case of vector valued functions. This is pointed out in next remark.

\begin{Remark}
Let $f\in C_{\mathrm{loc}}^{\frac{1}{2}+\gamma ,1}\left( \mathbb{R}_+ \times
\mathbb{R}^{n};\mathbb{R}^{p}\right) $\textrm{, }$g\in C_{\mathrm{loc}}^{%
\frac{1}{2}+\gamma ,1}\left( \mathbb{R}_+\times \mathbb{R}^{m};\mathbb{R}%
^{q}\right) $\textrm{, }$X=\left( X^{1},...,X^{n}\right) ^{\ast }$\textrm{, }%
$Y=\left( Y^{1},...,Y^{m}\right) ^{\ast }$\textrm{\ such that }$\left(
X,Y\right) $\textrm{\ has all its mutual covariations. Let }$V^{1},V^{2}$%
\textrm{\ be bounded variation processes. Then}
\begin{equation*}
\left[ f\left( V^{1},X\right) ,g\left( V^{2},Y\right) \right]
_{s}=\int_{0}^{s}\partial _{x}f\left( V^{1},X\right) d\left[ X,Y^{\ast }%
\right] \partial _{x}g\left( V^{2},Y\right) ^{\ast }
\end{equation*}
\hfill \textrm{for every }$s\ge 0.$\hfill
\hbox{\hskip 6pt\vrule width6pt height7pt
depth1pt  \hskip1pt}\bigskip
\end{Remark}

\smallskip One refined result is the vector It\^o formula whose proof
follows similarly as in \cite{RV2,FR}, where the involved stochastic
integrals were scalar.

\smallskip

\begin{Proposition}
\label{pr:ITOregular} Let $f\in C^{1,2}\left( \mathcal{T}_{t}\times \mathbb{R%
}^{n}\right) $. Let $X=\left( X^{1},...,X^{n}\right) ^{\ast }$, having all
its mutual covariations, $V$ be a bounded variation process indexed by $%
\mathcal{T}_{t}$. Then, for $0 \le s \le t$,
\begin{eqnarray*}
f\left( V_{s},X_{s}\right) &=&f\left( V_{t},X_{t}\right)
+\int_{t}^{s}\partial _{x}f\left( V_{r},X_{r}\right) d^{-}X_{r} \\
&&+\int_{t}^{s}\partial _{v}f\left( V_{r},X_{r}\right) dV_{r} \\
&&+\frac{1}{2}\int_{t}^{s}\partial _{xx}f\left( V_{r},X_{r}\right) d\left[
X,X^{\ast }\right] .
\end{eqnarray*}
\end{Proposition}

\begin{Remark}
\textrm{From the above statement it follows that, in particular, the
integral }$\int_{t}^{s}\partial _{x}f\left( V_{r},X_{r}\right) d^{-}X_{r}$%
\textrm{\ automatically exists.}\hfill
\hbox{\hskip 6pt\vrule width6pt height7pt
depth1pt  \hskip1pt}\bigskip
\end{Remark}

\begin{Remark}
\textrm{Let }$W=\left( W^{1},...,W^{n}\right) ^{\ast }$\textrm{\ be an }$%
\left( \mathcal{F}_{s}\right) $\textrm{-Brownian motion. Then }$\left[
W,W^{\ast }\right] _{s}=\left( \delta _{i,j}\right) s$\textrm{.\hfill
\hbox{\hskip 6pt\vrule width6pt height7pt
depth1pt  \hskip1pt}\bigskip }
\end{Remark}

\section{Fukushima-Dirichlet decomposition\label{ITO}}

\subsection{Definitions and remarks}

Throughout all this section we fix, as above, $0\leq t\leq T\leq +\infty $
and set $\mathcal{T}_{t}=\left[ t,T\right] \cap \mathbb{R}$. Recall that all
processes we consider are continuous except when explicitely stated.

\begin{Definition}
\label{df:Dirstrong} A real process $D$ is called an $\left( \mathcal{F}%
_{s}\right) $-Dirichlet process in $\mathcal{T}_{t}$ if it is $\left(
\mathcal{F}_{s}\right) $-adapted and can be written as
\begin{equation}
D=M+A,  \label{eq:dec1}
\end{equation}
where

\begin{itemize}
\item[(i)]  $M$ is an $\left( \mathcal{F}_{s}\right) $- local martingale,

\item[(ii)]  $A$ is a zero quadratic variation process such that (for
convenience) $A_{0}=0$.
\end{itemize}

A vector $D=\left( D^{1},...,D^{n}\right) $ is said to be Dirichlet if it
has all its mutual covariations and every $D^{i}$ is Dirichlet.
\end{Definition}

\begin{Remark}
\textrm{An }$\left( \mathcal{F}_{s}\right) $\textrm{-semimartingale is an }$%
\left( \mathcal{F}_{s}\right) $\textrm{-Dirichlet process.}\hfill
\hbox{\hskip 6pt\vrule width6pt height7pt
depth1pt  \hskip1pt}\bigskip
\end{Remark}

\begin{Remark}
\textrm{The decomposition (\ref{eq:dec1}) is unique, see for instance \cite
{RVW}.}\hfill \hbox{\hskip 6pt\vrule width6pt height7pt
depth1pt  \hskip1pt}\bigskip
\end{Remark}

The concept of Dirichlet process can be weakened for our purposes. We will
make use of an extension of such processes, called \textit{weak Dirichlet
processes}, introduced parallelly in \cite{ERbis} and implicitely in \cite{ER}.
 Recent developments  about the subject appear in \cite{cjms, RVSem, cr}.

 Weak Dirichlet processes are not Dirichlet processes but
they preserve a sort of orthogonal decomposition. In  all the papers
mentioned above however one deals with one-dimensional weak
Dirichlet processes while here
we treat the multidimensional case.%

\begin{Definition}
\label{df:Dirweak}A real process $D$ is called $\left( \mathcal{F}%
_{s}\right) $-\textbf{weak Dirichlet} process in $\mathcal{T}_{t}$ if it
can, be written as
\begin{equation}
D=M+A,  \label{eq:dec2}
\end{equation}
where

\begin{itemize}
\item[(i)]  $M$ is an $\left( \mathcal{F}_{s}\right) $-local martingale,

\item[(ii)]  $A$ is a process such that $\left[ A,N\right] =0$ for every $%
\left( \mathcal{F}_{s}\right) $-continuous local martingale $N$. (For
convenience $A_{0}=0$)
\end{itemize}

$A$ will be said \textbf{weak zero energy} process.

\end{Definition}

\begin{Remark}
The decomposition (\ref{eq:dec2}) is unique. In fact, let
\begin{equation*}
D=M^{1}+A^{1}=M^{2}+A^{2},
\end{equation*}
where $M^{1},M^{2},A^{1},A^{2}$ fulfill properties (i) and (ii) of
Definition \ref{df:Dirweak}. Then we have $M+A=0$ where $M=M^{1}-M^{2}$ and $%
A=A^{1}-A^{2}$ is such that $\left[ A,N\right] =0$ for every $\left(
\mathcal{F}_{s}\right) $- local martingale $N$.

It is now enough to evaluate the covariation of both members with $M$ to get
$\left[ M,M\right] =0$. Since $A_{0}=0$\ then $M_{0}=0$\ and consequently $%
M\equiv 0$. \hfill\qed
\end{Remark}

\begin{Example}
A simple example of weak Dirichlet process is given by a process $Z$ which
is independent of $\mathcal{F}$, for instance a deterministic one! Clearly
if $Z$ is not at least a finite quadratic variation process, it cannot be
Dirichlet. However it is possible to show that $[Z,N]=0$ for any local $%
\mathcal{F}$-martingale. In fact
\begin{eqnarray*}
\int_{0}^{s}(Z_{r+\varepsilon }-Z_{r})(N_{r+\varepsilon }-N_{r})dr
&=&\int_{0}^{s}dr(Z_{r+\varepsilon }-Z_{r})\int_{r}^{r+\varepsilon }\frac{1}{%
\varepsilon }dN_{\lambda } \\
&=&\int_{0}^{s}\frac{dN_{\lambda }}{\varepsilon }\int_{(\lambda -\varepsilon
)\vee 0}^{\lambda }(Z_{r+\varepsilon }-Z_{r})dr \\
&=&\int_{0}^{s}\frac{dN_{\lambda }}{\varepsilon }\int_{(\lambda -\varepsilon
)\vee 0}^{\lambda }Z_{r+\varepsilon }dr-\int_{0}^{s}\frac{dN_{\lambda }}{%
\varepsilon }\int_{(\lambda -\varepsilon )\vee 0}^{\lambda }Z_{r}dr.
\end{eqnarray*}
Previous expression converges u.c.p. to zero since the two last terms
converge u.c.p. to the It\^{o} integral $\int_{0}^{s}ZdN$ since $N$ is also
a local martingale with respect to the filtration $\mathcal{F}$ enlarged
with $Z$.
\end{Example}

\begin{Remark}
\label{concolution} \cite{ERbis} provides an example of weak Dirichlet
process coming from convolutions of local martingales. If for every $s>0$, $%
(G(s,\cdot ))$ is a continuous random field such that $(G(s,\cdot ))$ is $(%
\mathcal{F}_{r})$-progressively measurable and $M$ is an $(\mathcal{F}_{r})$%
-local martingale, then $X_{s}=\int_{t}^{s}G(s,r)dM_{r}$ defines a
weak Dirichlet process.\hfill\qed
\end{Remark}
\begin{Remark} \label{RBrownian}
\cite{cr} made the following observation.
If the underlying filtration  $\left( \mathcal{F}_{s}\right) $ is
 the natural filtration associated with a Brownian motion $(W_t)$
then  condition (ii) in Definition \ref{df:Dirweak} can be replaced with
\begin{description}
\item{(ii')}
$A$ is a process that $ [A,W] =0$. \hfill\qed
 \end{description}
\end{Remark}

\begin{Definition}
\label{weakvector} A vector $D=\left( D^{1},...,D^{n}\right) $ is said to be
a $\left( \mathcal{F}_{s}\right) $-\textbf{weak Dirichlet} process if every $%
D^{i}$ is $\left( \mathcal{F}_{s}\right) $- weak Dirichlet process. A vector
$A=\left( A^{1},...,A^{n}\right) $ is said to be a $\left( \mathcal{F}%
_{s}\right) $-\textbf{weak zero energy} process if every $A^{i}$ is $\left(
\mathcal{F}_{s}\right) $- weak zero energy process.
\end{Definition}

The aim now is to study what happens to a Dirichlet process after a $C^{0,1}$
type transformation. It is well known, see for instance \cite{RVW}, that a $%
C^{1}$ function of a finite quadratic variation process (respectively
Dirichlet process) is a finite quadratic variation (respectively Dirichlet)
process. Here, motivated by applications to optimal control, see the
introduction, section \ref{INTRO}, we look at the (possibly) inhomogeneous
case showing two different results: the first result (Proposition \ref
{pr:FDdec2} with Corollary \ref{Co:FDdec2}) states that a function $%
C^{0,1}\left( \mathbb{R_{+}}\times \mathbb{R}^{n}\right) $ of a weak
Dirichlet (vector) process having all its mutual covariations is again a
weak Dirichlet process; the second (Proposition \ref{pr:FDdec1}) gives
stability for Dirichlet processes in the inhomogeneous case for functions in
$C_{\mathrm{loc}}^{\frac{1}{2}+\gamma ,1}\left( \mathcal{T}_{t}\times
\mathbb{R}^{n}\right) $.

The first result comes from the need of treating optimal control
problems where the state process is a semimartingale that solves a
classical SDE's, which is the case we treat in the companion paper
\cite{GR2}.

\subsection{The decomposition for $C^{0,1}$ functions}

We now go on with a result concerning weak Dirichlet processes. Suppose $%
\left( D_{s}\right) _{s\in \mathcal{T}_{t}}$ to be an $\left( \mathcal{F}%
_{s}\right) $-Dirichlet process with decomposition (\ref{eq:dec1}) where $A$
is a zero quadratic variation process. Given a $C^{0,1}$ function $u$ of $D$%
, we cannot expect that $Z=u(\cdot ,D)$ is a Dirichlet process. However one
can hope that it is at least a weak Dirichlet process. Indeed this result is
true even if $D$ is a weak Dirichlet process with finite quadratic variation.

\begin{Proposition}
\label{pr:FDdec2} Suppose $\left( D_{s}\right) _{s\geq 0}$ to be an $\left(
\mathcal{F}_{s}\right) $- weak Dirichlet (vector) process having all its
mutual covariations. For every $u\in C^{0,1}\left( \mathbb{R}_{+}\times
\mathbb{R}^{n}\right) $ we have, for $s\geq t$,
\begin{equation}
u\left( s,D_{s}\right) =u\left( t,D_{t}\right) +\int_{t}^{s}\partial
_{x}u\left( r,D_{r}\right) dM_{r}+\mathcal{B}^{D}\left( u\right) _{s}-%
\mathcal{B}^{D}\left( u\right) _{t},  \label{eq:decweak}
\end{equation}
\smallskip where $\mathcal{B}^{D}:C^{0,1}\left( \mathbb{R}_{+}\times \mathbb{%
R}^{n}\right) \rightarrow \mathcal{C}_{\mathcal{F}}\left( \mathbb{R}%
_{+}\times \Omega ;\mathbb{R}^{n}\right) $ is a linear map having the
following properties:

\begin{enumerate}
\item[a)]  $\mathcal{B}^{D}$ is continuous;

\item[b)]  if $u\in C^{1,2}\left( \mathbb{R}_{+}\times \mathbb{R}^{n}\right)
$ then
\begin{equation*}
\mathcal{B}^{D}\left( u\right) _{s}=\int_{0}^{s}\partial _{s}u\left(
r,D_{r}\right) dr+\int_{0}^{s}\partial _{xx}u\left( r,D_{r}\right)
d\left[ M,M\right] _{r}+\int_{0}^{s}\partial _{x}u\left(
r,D_{r}\right) d^{-}A_{r};
\end{equation*}

\item[c)]  if $u\in C^{0,1}\left( \mathbb{R}_{+}\times \mathbb{R}^{n}\right)
$ then $\left( (\mathcal{B} ^{D}\left( u\right) _{s} \right) $
is a weak zero energy process.
\end{enumerate}
\end{Proposition}

\begin{Corollary}
\label{Co:FDdec2} Suppose $\left( D_{t}\right) _{t\geq 0}$ to be an $\left(
\mathcal{F}_{t}\right) $- weak Dirichlet (vector) process having all its
mutual covariations. For every $u\in C^{0,1}\left( \mathbb{R}_{+}\times
\mathbb{R}^{n}\right) ,$ $(u(t,D_{t})$ is a $\left( \mathcal{F}_{t}\right) $%
-weak Dirichlet process with martingale part $\tilde{M}_{t}=\int_{0}^{t}%
\partial _{x}u(s,D_{s})dM_{s}$.
\end{Corollary}

\begin{Remark}
\label{rm:FDdec2} Given a bounded stopping time $\tau $ with values
in $\mathcal{T}_{t}$, it is easy to see that decomposition
(\ref{eq:decweak}) still holds for the stopped process $D^{\tau }$.
In fact given an $\left( \mathcal{F}_{s}\right) $- martingale $N$,
and an $(\mathcal{F}_{s})$- weak zero energy process $A$, we have
$[N,A^{\tau }]=[N,A]^{\tau }=0$.\hfill\qed
\end{Remark}

\textbf{Proof} (of the Proposition).
Without restriction of generality we will set $t=0$.

Property a) follows simply by writing
\begin{equation*}
\mathcal{B}^{D}\left( u\right) _{s} = u\left( s,D_{s}\right) -u\left(
0,D_{0}\right) -\int_{0}^{s}\partial _{x}u\left( r,D_{r}\right) dM_{r}
\end{equation*}
and observing that the process defined on the right hand side has the
required continuity property.

Property b) follows from Proposition \ref{pr:ITOregular} applied reversely.
Indeed, given $u\in C^{1,2}\left( \mathbb{R}_+\times \mathbb{R}^{n}\right) $%
, Proposition \ref{pr:ITOregular} can be applied. In particular
\begin{equation*}
\int_{0}^{s}\partial _{x}u\left( r,D_{r}\right) d^{-}D_{r}
\end{equation*}
exists; this implies that also
\begin{equation*}
\int_{0}^{s}\partial _{x}u\left( r,D_{r}\right) d^{-}A_{r}
\end{equation*}
exists since
\begin{equation*}
\int_{0}^{s}\partial _{x}u\left( r,D_{r}\right)
d^{-}M_{r}=\int_{0}^{s}\partial _{x}u\left( r,D_{r}\right) dM_{r}
\end{equation*}
is the classical It\^{o} integral.

It remains to prove point c)

\begin{equation*}
\left[ u\left( \cdot ,D\right) -\int_{0}^{\cdot }\partial _{x}u\left(
r,D_{r}\right) dM_{r},N\right] =0
\end{equation*}
for every one dimensional $\left( \mathcal{F}_{r}\right) $- local martingale
$N$.

For simplicity of notations, we will suppose that $D$ is one-dimensional.
Therefore $D$ will be therefore a finite quadratic variation proces. Since
the covariation of semimartingales coincides with the classical covariation
\begin{equation*}
\left[ \int_{0}^{\cdot }\partial _{x}u\left( r,D_{r}\right) dM_{r},N\right]
=\int_{0}^{\cdot }\partial _{x}u\left( r,D_{r}\right) d\left[ M,N\right]
_{r},
\end{equation*}
it remains to check that, for every $s\in \left[ 0,T\right] $,
\begin{equation*}
\left[ u\left( \cdot ,D\right) ,N\right] _{s}=\int_{0}^{s}\partial
_{x}u\left( r,D_{r}\right) d\left[ M,N\right] _{r}.
\end{equation*}
For this, we have to evaluate the $u.c.p.$ limit of
\begin{equation*}
\int_{0}^{s}\left[ u\left( r+\varepsilon ,D_{r+\varepsilon }\right) -u\left(
r,D_{r}\right) \right] \frac{N_{r+\varepsilon }-N_{r}}{\varepsilon }dr
\end{equation*}
in probability. This can be written as the sum of two terms:
\begin{equation*}
I_{1}\left( s,\varepsilon \right) =\int_{0}^{s}\left( u\left( r+\varepsilon
,D_{r+\varepsilon }\right) -u\left( r+\varepsilon ,D_{r}\right) \right)
\frac{N_{r+\varepsilon }-N_{r}}{\varepsilon }dr,
\end{equation*}
\begin{equation*}
I_{2}\left( s,\varepsilon \right) =\int_{0}^{s}\left( u\left( r+\varepsilon
,D_{r}\right) -u\left( r,D_{r}\right) \right) \frac{N_{r+\varepsilon }-N_{r}%
}{\varepsilon }dr.
\end{equation*}
First we prove that $I_{1}\left( s,\varepsilon \right) $ goes to $%
\int_{0}^{t}\partial _{x}u\left( r,D_{r}\right) d\left[ M,N\right] _{r}$. In
fact
\begin{eqnarray}
I_{1}\left( s,\varepsilon \right) &=&\int_{0}^{s}\left( u\left(
r+\varepsilon ,D_{r+\varepsilon }\right) -u\left( r+\varepsilon
,D_{r}\right) \right) \frac{N_{r+\varepsilon }-N_{r}}{\varepsilon }dr  \notag
\\
&=&\int_{0}^{s}\partial _{x}u\left( r+\varepsilon ,D_{r}\right) \left(
D_{r+\varepsilon }-D_{r}\right) \frac{N_{r+\varepsilon }-N_{r}}{\varepsilon }%
dr+R_{1}\left( s,\varepsilon \right) ,  \label{eq:splitI1}
\end{eqnarray}
where $R_{1}\left( s,\varepsilon \right) \rightarrow 0$ $u.c.p$\textit{.} as
$\varepsilon \rightarrow 0$. Indeed
\begin{eqnarray*}
R_{1}\left( s,\varepsilon \right) &=&\int_{0}^{s}\left[ \int_{0}^{1}\left[
\partial _{x}u\left( r+\varepsilon ,D_{r}+\lambda \left( D_{r+\varepsilon
}-D_{r}\right) \right) -\partial _{x}u\left( r+\varepsilon ,D_{r}\right) %
\right] d\lambda \right. \\
&&\left. \frac{N_{r+\varepsilon }-N_{r}}{\varepsilon }\left(
D_{r+\varepsilon }-D_{r}\right) \right] dr
\end{eqnarray*}
and the claim follows by the continuity of $\partial _{x}u$ and from the
estimate
\begin{eqnarray}
&&\frac{1}{\varepsilon }\int_{0}^{T}\left( N_{r+\varepsilon }-N_{r}\right)
\left( D_{r+\varepsilon }-D_{r}\right) dr  \notag \\
&\leq &\left[ \frac{1}{\varepsilon }\int_{0}^{T}\left( N_{r+\varepsilon
}-N_{r}\right) ^{2}dr\cdot \frac{1}{\varepsilon }\int_{0}^{T}\left(
D_{r+\varepsilon }-D_{r}\right) ^{2}dr\right] ^{\frac{1}{2}}\overset{%
\varepsilon \rightarrow 0}{\longrightarrow }\left( \left[ N\right] \cdot %
\left[ D\right] \right) ^{\frac{1}{2}}.  \label{eq:convNDvarquad}
\end{eqnarray}
On the other hand the first term in (\ref{eq:splitI1}) can be rewritten as
\begin{equation}
\int_{0}^{s}\partial _{x}u\left( r,D_{r}\right) \left( D_{r+\varepsilon
}-D_{r}\right) \frac{N_{r+\varepsilon }-N_{r}}{\varepsilon }dr+R_{2}\left(
s,\varepsilon \right) ,  \label{eq:A4}
\end{equation}
where $R_{2}\left( s,\varepsilon \right) \rightarrow 0$ $u.c.p.$ arguing as
for $R_{1}\left( s,\varepsilon \right) $. The integral in (\ref{eq:A4}) goes
then $u.c.p.$ to $\int_{0}^{s}\partial _{x}u\left( r,D_{r}\right) d\left[ M,N%
\right] _{r}$, since by (\ref{eq:convNDvarquad}) the measures $\frac{\left(
N_{r+\varepsilon }-N_{r}\right) \left( D_{r+\varepsilon }-D_{r}\right) }{%
\varepsilon }dr$ weakly converge to $d\left[ N,D\right] $ as $\varepsilon
\rightarrow 0$.

It remains to show that $I_{2}\left( s,\varepsilon \right) \rightarrow 0$ $%
u.c.p.$ for every $s\in \left[ 0,T\right] $ as $\varepsilon \rightarrow 0$.
By using suitable localization theorems (e.g. as usually done for instance
in \cite{RY}, section IV.1), it is enough to suppose $u$ to be with compact
support and $N$ to be a square integrable martingale. Then we evaluate

\begin{equation}
\mathbb{E}\left( \sup_{0\leq s\leq T}|I_{2}(s,\varepsilon )|^{2}\right) .
\label{I2}
\end{equation}
Now we have, exchanging integrals
\begin{eqnarray*}
I_{2}(s,\varepsilon ) &=&\int_{0}^{s}[u(r+\varepsilon ,D_{r})-u(r,D_{r})]%
\frac{N_{r+\varepsilon }-N_{r}}{\varepsilon }dr \\
&=&\int_{0}^{s}[u(r+\varepsilon ,D_{r})-u(r,D_{r})]dr\int_{r}^{r+\varepsilon
}\frac{1}{\varepsilon }dN_{\lambda } \\
&=&\int_{0}^{s}\frac{dN_{\lambda }}{\varepsilon }\int_{(\lambda -\varepsilon
)\vee 0}^{\lambda }[u(r+\varepsilon ,D_{r})-u(r,D_{r})]dr.
\end{eqnarray*}
Doob inequality implies that (\ref{I2}) is smaller than
\begin{equation*}
4\mathbb{E}\left\{ \int_{0}^{T}d[N]_{\lambda }\left( \frac{1}{\varepsilon }%
\int_{(\lambda -\varepsilon )\vee 0}^{\lambda }[u(r+\varepsilon
,D_{r})-u(r,D_{r})]dr\right) ^{2}\right\} .
\end{equation*}
The fact that $u$ is uniformly continuous on compact sets and the Lebesgue
dominated convergence theorem imply the result.\hfill
\hbox{\hskip 6pt\vrule width6pt height7pt
depth1pt  \hskip1pt}\bigskip

A significant bracket evaluation, in the spirit of Proposition \ref
{pr:bracket}, but for $C^{0,1}-$ functions of semimartingales is the
following. For simplicity we formulate the one-dimensional case, even if it
extends to the multidimensional case.

\begin{Corollary}
\label{cor:bracket} Let $S$ be a $\left( \mathcal{F}_{s}\right)$-
semimartingale in $\mathbb{R}_+$, $u\in C^{0,1}( \mathcal{T}_{t}\times
\mathbb{R}) $. Then

\begin{equation*}
\lbrack u(\cdot ,S),S]_{t}=\int_{0}^{t}\partial _{x}u(s,S_{s})d[S]_{s}.
\end{equation*}
\end{Corollary}

\textbf{Proof.} Let $S= M + V$ be the decomposition of $S$ with $M$ being a
local martingale and $V$ a finite variation process with $V_0 = 0.$ Then
\begin{equation*}
u(t,S_t) = u(0,S_0) + \int_0^t \partial_x u(s,S_s) dM_s + \tilde A_t,
\end{equation*}
where $\tilde A$ is a weak zero energy process. In particular, a classical
localization argument shows that $[\tilde A, M] = 0$. On the other hand,
obviously $[\tilde A, V] = 0$; consequently, by linearity and since the
covariation of local martingales is the classical convolution, the result
follows. \hfill
\hbox{\hskip 6pt\vrule width6pt height7pt
depth1pt  \hskip1pt}\bigskip

\subsection{The decomposition for $C^{\frac{1}{2}+\protect\gamma ,1}$
functions\label{DECHOLDER}}

If, in Proposition \ref{pr:FDdec2}, $D$ is a Dirichlet process and $u$ is of
class $C^{\frac{1}{2}+\gamma ,1}$, $\gamma >0$, then the results of
Proposition \ref{pr:FDdec2} and Corollary \ref{Co:FDdec2} can be better
precised. In fact it is possible to show that Dirichlet processes are stable
through $C^{\frac{1}{2}+\gamma ,1}$, $\gamma >0$ transformations.

\begin{Proposition}
\label{pr:FDdec1}Let $\left( D_{s}\right) _{s\geq 0}$ be an $\left( \mathcal{%
F}_{s}\right) $-Dirichlet process with decomposition (\ref{eq:dec1}). The
statement of Proposition \ref{pr:FDdec2} holds with
\begin{equation*}
\mathcal{B}^{D}:C^{\frac{1}{2}+\gamma ,1}\left( \mathbb{R}_{+}\times \mathbb{%
R}^{n}\right) \rightarrow \mathcal{C}_{\mathcal{F}}\left( \mathbb{R}%
_{+}\times \Omega ;\mathbb{R}^{n}\right)
\end{equation*}
fulfilling properties a), b) and

{c)} if $u\in C^{\frac{1}{2}+\gamma ,1}\left( \mathbb{R}_{+}\times \mathbb{R}%
^{n}\right) $ then $\left( \mathcal{B}^{D}\left( u\right) _{s}\right) $ is a
zero quadratic variation process.
\end{Proposition}

\textbf{Proof. } Points a) and b) follow similarly as for the $C^{0,1}$
decomposition.

In order to establish Property c) we proceed using the bilinearity of the
covariation. We will in fact show that $\mathcal{B}^{D}\left( u \right)$
is a zero quadratic variation process. We operate with the bilinearity of
the covariation process and we evaluate

\begin{enumerate}
\item[(i)]  $\left[ u\left( \cdot ,D\right) ,u\left( \cdot ,D\right) \right]
,$

\item[(ii)]  $\left[ u\left( \cdot ,D\right) ,\int_{0}^{\cdot }\partial
_{x}u\left( r,D_{r}\right) dM_{r}\right] ,$

\item[(iii)]  $\left[ \int_{0}^{\cdot }\partial _{x}u\left( r,D_{r}\right)
dM_{r},\int_{0}^{\cdot }\partial _{x}u\left( r,D_{r}\right) dM_{r}\right] ,$
\end{enumerate}

as follows.

\begin{enumerate}
\item[(i)]  We apply Proposition \ref{pr:bracket} to get that
\begin{equation*}
\left[ u\left( \cdot ,D\right) ,u\left( \cdot ,D\right) \right]
_{s}=\int_{0}^{s}\partial _{x}u\left( s,D_{r}\right) d\left[ D,D^{\ast }%
\right] _{r}\partial _{x}u\left( r,D_{r}\right) ^{\ast }.
\end{equation*}

\item[(ii)]  Setting $N_{t}=\int_{0}^{s}\partial _{x}u\left( r,D_{r}\right)
dM_{r}$, Remark \ref{rm:bracketnullo} implies that $\left( N,D\right) $ has
all its mutual brackets; therefore again Proposition \ref{pr:bracket}
implies that
\begin{equation*}
\left[ u\left( \cdot ,D\right) ,N^{\ast }\right] _{s}=\int_{0}^{s}\partial
_{x}u\left( r,D_{r}\right) d\left[ D,N^{\ast }\right] _{r}.
\end{equation*}
On the other hand, by Remark \ref{rm:bracketnullo}
\begin{equation*}
\left[ D,N^{\ast }\right] _{t}=\left[ M,N^{\ast }\right] _{s}=\int_{0}^{s}%
\partial _{x}u\left( r,D_{r}\right) d\left[ M,M^{\ast }\right] _{r}\partial
_{x}u\left( r,D_{r}\right) ^{\ast }.
\end{equation*}

\item[(iii)]  The fact that the covariation of semimartingales coincides
with the classical covariation gives
\begin{eqnarray*}
&&\left[ \int_{0}^{\cdot }\partial _{x}u\left( r,D_{r}\right)
dM_{r},\int_{0}^{\cdot }\partial _{x}u\left( r,D_{r}\right) dM_{r}\right]
_{s} \\
&=&\int_{0}^{s}\partial _{x}u\left( r,D_{r}\right) d\left[ M,M^{\ast }\right]
_{r}\partial _{x}u\left( r,D_{r}\right) ^{\ast }.
\end{eqnarray*}
Finally by Remark \ref{rm:bracketnullo} and the decomposition we get that
\begin{equation*}
\left[ D,D^{\ast }\right] =\left[ M,M^{\ast }\right] .
\end{equation*}
\end{enumerate}

The bilinearity of the covariation allows now to conclude.\hfill
\hbox{\hskip 6pt\vrule width6pt height7pt
depth1pt  \hskip1pt}\bigskip

\smallskip

\section{Representation of operator $\mathcal{B}$ when $u$ solves a suitable
PDE.\label{REPRESENTATION}}

Here we want to develop the connection between suitable (deterministic)
linear differential operators and our (stochastic) operators $\mathcal{B}$
introduced in the previous section. This connection is well known and
obvious when $u$ is the $C^{1,2}$ solution of a second order PDE and $D$ is
a diffusion process.
Our aim is to extend the validity of such representation when $u$ is only a $%
C^{0,1}$ solution (in a suitable sense that we will define below) and $D$ is
a weak Dirichlet process of a suitable kind (see below). This will be used
as a key tool in the applications to optimal control.

\subsection{Strong solutions of parabolic PDE's}

Let $0<T<+\infty $, consider two continuous functions
\begin{equation*}
b:\left[ 0,T\right] \times \mathbb{R}^{n}\rightarrow \mathbb{R}^{n},\qquad
\sigma :\left[ 0,T\right] \times \mathbb{R}^{n}\rightarrow \mathcal{L}\left(
\mathbb{R}^{m},\mathbb{R}^{n}\right)
\end{equation*}
and the linear parabolic operator
\begin{equation*}
\mathcal{L}_{0}:D\left( \mathcal{L}_{0}\right) \subseteq C^{0}\left( \left[
0,T\right] \times \mathbb{R}^{n}\right) \longrightarrow C^{0}\left( \left[
0,T\right] \times \mathbb{R}^{n}\right) ,
\end{equation*}
\begin{equation*}
D\left( \mathcal{L}_{0}\right) =C^{1,2}\left( \left[ 0,T\right] \times
\mathbb{R}^{n}\right) ,
\end{equation*}
\begin{equation*}
\mathcal{L}_{0}u\left( t,x\right) =\partial _{t}u\left( t,x\right)
+\left\langle b\left( t,x\right) ,\partial _{x}u\left( t,x\right)
\right\rangle +\frac{1}{2}\text{Tr }\left[ \sigma ^{\ast }\left( t,x\right)
\partial _{xx}u\left( t,x\right) \sigma \left( t,x\right) \right] .
\end{equation*}
Defining
\begin{equation*}
L_{0}\left( t\right) u\left( t,x\right) =\left\langle b\left( t,x\right)
,\partial _{x}u\left( t,x\right) \right\rangle +\frac{1}{2}\text{Tr }\left[
\sigma ^{\ast }\left( t,x\right) \partial _{xx}u\left( t,x\right) \sigma
\left( t,x\right) \right] ,
\end{equation*}
we can write
\begin{equation*}
\mathcal{L}_{0}u\left( t,x\right) =\partial _{t}u\left( t,x\right)
+L_{0}\left( t\right) u\left( t,x\right) .
\end{equation*}
Recall that an operator $\mathcal{M}:D\left( \mathcal{M}\right) \subseteq
F\rightarrow G$ ($F,G$ suitable Fr\'{e}chet spaces) is closable if, given
any sequence $\left( u_{n}\right) _{n\in \mathbb{N}}\subseteq D\left(
\mathcal{M}\right) $, we have
\begin{equation*}
\left.
\begin{array}{l}
u_{n}\longrightarrow 0,\quad \quad \text{in }F \\
\mathcal{M}u_{n}\longrightarrow \eta \quad \text{in }G
\end{array}
\right\} \Longrightarrow \eta =0.
\end{equation*}

When $\mathcal{L}_{0}$ is closable we denote by $\mathcal{L}$ its closure
and recall that
\begin{equation*}
u\in D\left( \mathcal{L}\right) \Longleftrightarrow \exists \left(
u_{n}\right) \subset C^{1,2}\left( \left[ 0,T\right] \times \mathbb{R}%
^{n}\right) :\left\{
\begin{array}{l}
u_{n}\longrightarrow u \\
\mathcal{L}_{0}u_{n}\longrightarrow \mathcal{L}u
\end{array}
\right. \text{in }C^{0}\left( \left[ 0,T\right] \times \mathbb{R}^{n}\right)
.
\end{equation*}

Now, given $\phi \in C^{0}\left( \mathbb{R}^{n}\right) $ and $h\in
C^{0}\left( \left[ 0,T\right] \times \mathbb{R}^{n}\right) $ we consider the
inhomogenous backward parabolic problem
\begin{equation}
\left\{
\begin{array}{l}
\partial _{t}u\left( t,x\right) +\left\langle b\left( t,x\right) ,\partial
_{x}u\left( t,x\right) \right\rangle +\frac{1}{2}\text{Tr }\left[ \sigma
^{\ast }\left( t,x\right) \partial _{xx}u\left( t,x\right) \sigma \left(
t,x\right) \right] =h\left( t,x\right) , \\
\qquad \qquad \qquad \qquad \qquad \qquad \qquad \qquad \qquad \qquad t\in
\left[ 0,T\right] \text{,\qquad }x\in \mathbb{R}^{n}, \\
u\left( T,x\right) =\phi \left( x\right) ,\qquad \text{\qquad }x\in \mathbb{R%
}^{n},
\end{array}
\right.  \label{eq:CPlinear}
\end{equation}
that can be rewritten as
\begin{equation*}
\partial _{t}u\left( t,x\right) +L_{0}\left( t\right) u\left( t,x\right)
=h\left( t,x\right) ,\qquad u\left( T,x\right) =\phi \left( x\right) ,
\end{equation*}
or
\begin{equation*}
\mathcal{L}_{0}u\left( s,x\right) =h\left( s,x\right) ,\qquad u\left(
T,x\right) =\phi \left( x\right) .
\end{equation*}

\begin{Definition}
\label{df:solstrict}We say that $u\in C^{0}\left( \left[ 0,T\right] \times
\mathbb{R}^{n}\right) $ is a strict solution to the backward Cauchy problem (%
\ref{eq:CPlinear}) if $u\in D\left( \mathcal{L}_{0}\right) $ and (\ref
{eq:CPlinear}) holds.
\end{Definition}

\begin{Definition}
\label{df:solstrong}We say that $u\in C^{0}\left( \left[ 0,T\right] \times
\mathbb{R}^{n}\right) $ is a strong solution to the backward Cauchy problem (%
\ref{eq:CPlinear}) if there exists a sequence $\left( u_{n}\right) \subset
D\left( \mathcal{L}_{0}\right) $ and two sequences $\left( \phi _{n}\right)
\subseteq C^{0}\left( \mathbb{R}^{n}\right) $, $\left( h_{n}\right)
\subseteq C^{0}\left( \left[ 0,T\right] \times \mathbb{R}^{n}\right) $, such
that

\begin{enumerate}
\item  For every $n\in \mathbb{N}$ $u_{n}$ is a strict solution of the
problem
\begin{equation*}
\mathcal{L}_{0}u_{n}\left( t,x\right) =h_{n}\left( t,x\right) ,\qquad
u_{n}\left( T,x\right) =\phi _{n}\left( x\right) .
\end{equation*}

\item  The following limits hold
\begin{equation*}
\begin{array}{l}
u_{n}\longrightarrow u\text{ in }C^{0}\left( \left[ 0,T\right] \times
\mathbb{R}^{n}\right) , \\
h_{n}\longrightarrow h\text{ in }C^{0}\left( \left[ 0,T\right] \times
\mathbb{R}^{n}\right) , \\
\phi _{n}\longrightarrow \phi \text{ in }C^{0}\left( \mathbb{R}^{n}\right) .
\end{array}
\end{equation*}
\end{enumerate}
\end{Definition}

\bigskip

\subsection{The representation result}

Let $u$ be a strong solution of class $C^{0,1}$ of (\ref{eq:CPlinear}) and $%
S $ be a weak Dirichlet process that can be written in the following form:

\begin{equation*}
S_{t}=S_{0}+\int_{0}^{t}\sigma (s,S_{s})dW_{s}+A_{t}
\end{equation*}
where $\sigma $ is as in the previous section and $A$ is a weak zero energy
process with finite quadratic variation.

We observe that for our applications to optimal control it would be
enough to take $S$ semimartingale. We deal with this more general
case to prepare the field for a forthcoming paper in which we
consider the optimal control of solutions of SDEs where the drift
$b$ is the derivative in space of a continuous function $\beta $,
therefore a Schwartz distribution. One would have in that case
$A_{t}=``\int_{0}^{t}\partial _{x}\beta (s,x)ds"$ in some specific
sense. Equations of that type, when there is no dependence in time
appear for instance in \cite{frw2}. Solutions are Dirichlet
processes in the time-homogeneous case and weak Dirichlet in the
general case.

We remark that the coefficient $\sigma $ \emph{must coincide} with the one
appearing in the second order term of the operator.

We state first a technical lemma whose proof is elementary.

\begin{Lemma}
\label{lm:newbyFR}Let $T<+\infty $. Let $f_{n},f:\left[ 0,T\right]
\rightarrow \mathbb{R}$, $n\in \mathbb{N}$, continuous such that $%
f_{n}\rightarrow f$ uniformly. For a fixed constant $K>0$, we define
\begin{equation*}
\tau _{n}=\inf \left\{ t\in \left[ 0,T\right] :\left| f_{n}\left( t\right)
\right| \geq K\right\} ,\qquad \tau =\inf \left\{ t\in \left[ 0,T\right]
:\left| f\left( t\right) \right| \geq K\right\}
\end{equation*}
with the convention that $\inf \emptyset =T$. Then
\begin{equation*}
\lim_{n\rightarrow +\infty }f_{n}^{\tau _{n}}(T)=f^{\tau }(T)
\end{equation*}
where $f^{\tau }$ (respectively $f_{n}^{\tau _{n}}$) is the stopped
function defined by $f^{\tau }\left( t\right) =f \left( \tau \wedge
t\right) $ (respectively $f_{n}^{\tau _{n}}\left( t\right)
=f_{n}\left( \tau _{n}\wedge t\right) $).
\end{Lemma}

>From the above lemma we get the following Proposition.

\begin{Proposition}
\label{pr:FRnew}The set $\mathcal{M}_{loc}$ of all $\left(
\mathcal{F}_{t}\right) -$continuous local martingales is a closed
subset of $\mathcal{C}_{\mathcal{F}}\left( \left[ 0,T\right] \times
\Omega ;\mathbb{R}\right) $ endowed with the $u.c.p.$ topology.
\end{Proposition}

\textbf{Proof. }Let $\left( M_{n}\left( t\right) ,\quad t\in \left[ 0,T%
\right] \right) _{n\in \mathbb{N}}$ be a sequence of local
continuous martingales converging $u.c.p.$ to a continuous process
$M\in \mathcal{C}_{\mathcal{F}}\left( \left[ 0,T\right] \times
\Omega ;\mathbb{R}\right) $. For $K>0$ we define the following
stopping times:
\begin{equation*}
\tau ^{n}=\inf \left\{ t\in \left[ 0,T\right] :\left| M_{n}\left( t\right)
\right| \geq K\right\}
\end{equation*}
\begin{equation*}
\tau =\inf \left\{ t\in \left[ 0,T\right] :\left| M\left( t\right) \right|
\geq K\right\}
\end{equation*}
with the convention that $\inf \emptyset =T$.

In order to conclude, it suffices to show that $M^{\tau }$ is a square
integrable martingale. Lemma \ref{lm:newbyFR} implies that
\begin{equation}
M_{n}^{\tau ^{n}}(T)\longrightarrow M^{\tau }(T)\qquad a.e.
\label{eq:convmartFR}
\end{equation}
Using (\ref{eq:convmartFR}) above and Lebesgue dominated convergence
theorem, we have
\begin{equation*}
M_{n}^{\tau ^{n}}\left( T\right) \longrightarrow M^{\tau }\left( T\right)
\qquad \text{in }L^{2}\left( \Omega \right) .
\end{equation*}
The fact that $M^{\tau }$ is a square integrable martingale follows then
from Proposition 5.23, Ch. 1 of \cite{KS}.\hfill
\hbox{\hskip 6pt\vrule width6pt height7pt
depth1pt  \hskip1pt}\bigskip

We are now able to state a useful representation result. Below we fix $t\in %
\left[ 0,T\right] $ and, since now $T<+\infty $, we have $\mathcal{T}_{t}=%
\left[ t,T\right] $.

\begin{Theorem}
\label{th:Itostrong}Let $T<+\infty $ and
\begin{equation*}
b:\left[ 0,T\right] \times \mathbb{R}^{n}\rightarrow \mathbb{R}^{n},\qquad
\mathrm{resp.}\ \sigma :\left[ 0,T\right] \times \mathbb{R}^{n}\rightarrow
\mathcal{L}\left( \mathbb{R}^{m},\mathbb{R}^{n}\right)
\end{equation*}
be continuous functions. Let $u\in C^{0,1}\left( \left[ 0,T\right] \times
\mathbb{R}^{n}\right) $ be a strong solution of the Cauchy problem (\ref
{eq:CPlinear}).

Fix $t\in \left[ 0,T\right] $, $x\in \mathbb{R}^{n}$ and let $(S_{s})_{s\in
\mathcal{T}_{t}}$ be a process of the form
\begin{equation*}
S_{s}=x+\int_{t}^{s}\sigma \left( r,S_{r}\right) dW_{r}+A_{t}-A_{s}
\end{equation*}
where $(A_{s})_{s\in \mathcal{T}_{t}}$ is an $(\mathcal{F}_{s})$-weak zero
energy process having all its mutual covariations.

Then, provided that the following assumption is verified for every $s\in
\mathcal{T}_{t}$: 
\begin{eqnarray}
\lim_{n\rightarrow +\infty }\int_{t}^{s}(\partial _{x}u_{n}(r,S_{r})
&-&\partial _{x}u(r,S_{r}))d^{-}A_{r}  \label{eq:hpnewconvderucp} \\
&-&\int_{t}^{s}\left\langle \partial _{x}u_{n}(r,S_{r})-\partial
_{x}u(r,S_{r}),b(r,S_{r})\right\rangle dr=0,\qquad u.c.p.,  \notag
\end{eqnarray}
we have
\begin{equation}
u(s,S_{s})=u(t,S_{t})+\int_{t}^{s}\partial _{x}u\left( r,S_{r}\right) \sigma
(r,S_{r})dW_{r}+\mathcal{B}^{S}\left( u\right) _{s}-\mathcal{B}^{S}\left(
u\right) _{t},  \label{Estrong}
\end{equation}
where, for $s\in \mathcal{T}_{t}$%
\begin{equation}
\mathcal{B}^{S}\left( u\right) _{s}=\int_{0}^{s}h\left( r,S_{r}\right)
dr+\int_{0}^{s}\partial _{x}u\left( r,S_{r}\right)
d^{-}A_{r}-\int_{0}^{s}\left\langle \partial _{x}u\left( r,S_{r}\right)
,b\left( r,S_{r}\right) \right\rangle dr.  \label{eq:reprmain}
\end{equation}
\end{Theorem}

\textbf{Proof.} We set $t=0$ for simplicity. The general case is obtained by
additivity of different integrals. Let $u_{n}\longrightarrow u$ in $%
C^{0}\left( \left[ 0,T\right] \times \mathbb{R}^{n}\right) $ be a sequence
such that $\mathcal{L}_{0}u_{n}=h_{n}\longrightarrow h$ in $C^{0}\left( %
\left[ 0,T\right] \times \mathbb{R}^{n}\right) $. By Proposition \ref
{pr:ITOregular}, we get
\begin{eqnarray*}
u_{n}\left( s,S_{s}\right) &=&u_{n}\left( 0,S_{0}\right) + \int_{0}^{s}%
\mathcal{L}_{0}u_{n}\left( r,S_{r}\right) dr - \int_{0}^{s}\left\langle
\partial _{x}u_{n}\left( r,S_{r}\right) , b\left( r,S_{r}\right)
\right\rangle dr \\
&+& \int_{0}^{s}\partial _{x}u_n \left( r,S_{r}\right) \sigma \left(
r,S_{r}\right) dW_{r} + \int_{0}^{s} \partial _{x}u_{n}\left( r,S_{r}
\right) d^- A_r.
\end{eqnarray*}

>From (\ref{eq:hpnewconvderucp}) we conclude that
\begin{eqnarray*}
M_{s}^{n} &=&u_{n}\left( s,S_{s}\right) -u_{n}\left( 0,S_{0}\right)
-\int_{0}^{s}\mathcal{L}_{0}u_{n}\left( r,S_{r}\right) dr \\
&+&\int_{0}^{s}\left\langle \partial _{x}u_{n}\left( r,S_{r}\right) ,b\left(
r,S_{r}\right) \right\rangle dr - \int_{0}^{s} \partial _{x} u_{n}\left(
r,S_{r} \right) d^- A_r
\end{eqnarray*}
converges $u.c.p.$ to
\begin{eqnarray*}
M_{s} &=&u\left( s,S_{s}\right) -u\left( 0,S_{0}\right) -\int_{0}^{s}h\left(
r,S_{r}\right) dr \\
&+&\int_{0}^{s}\left\langle \partial _{x}u\left( r,S_{r}\right) ,b \left(
r,S_{r}\right) \right\rangle dr \int_{0}^{s} \partial _{x}u \left( r,S_{r}
\right) d^- A_.
\end{eqnarray*}
Using Proposition \ref{pr:FRnew} above we get that $M$ is an $\left(
\mathcal{F}_{s}\right) $-local martingale. The result follows by Proposition
\ref{pr:FDdec2}, where $dM = \sigma(s,S_s) ds$ and $D = S,$ by
identification of the weak zero energy processes.\hfill
\hbox{\hskip 6pt\vrule width6pt height7pt
depth1pt  \hskip1pt}\bigskip

For our applications in \cite{GR2} we will need to consider a process $A$
which is of bounded variation (so that $S$ is solution of an SDE) but which is
non-Markovian.

\begin{Corollary}
\label{cor:ITOstrong}Let $T<+\infty $ and
\begin{equation*}
b_{1}:\Omega \times \left[ 0,T\right] \times \mathbb{R}^{n}\rightarrow
\mathbb{R}^{n},
\end{equation*}
be a continuous progressively measurable field (continuous in $(s,x)$) and
\begin{equation*}
b:\left[ 0,T\right] \times \mathbb{R}^{n}\rightarrow \mathbb{R}^{n},\qquad
\sigma :\left[ 0,T\right] \times \mathbb{R}^{n}\rightarrow \mathcal{L}\left(
\mathbb{R}^{m},\mathbb{R}^{n}\right)
\end{equation*}
be continuous functions. Let $u\in C^{0,1}\left( \left[ 0,T\right] \times
\mathbb{R}^{n}\right) $ be a strong solution of the Cauchy problem (\ref
{eq:CPlinear}).

Fix $t\in \left[ 0,T\right] $, $x\in \mathbb{R}^{n}$ and let $(S_{s})$ be a
solution to the SDE
\begin{equation*}
dS_{s}=b_{1}\left( s,S_{s}\right) ds+\sigma \left( s,S_{s}\right)
dW_{s};\qquad S_{t}=x.
\end{equation*}

Then, provided that the following assumption be verified for every $s\in
\mathcal{T}_{t}$%
\begin{equation}
\lim_{n\rightarrow +\infty }\int_{t}^{s}\left\langle \partial
_{x}u_{n}\left( r,S_{r}\right) -\partial _{x}u\left( r,S_{r}\right)
,b_{1}\left( r,S_{r}\right) -b\left( r,S_{r}\right) \right\rangle
dr=0,\qquad u.c.p.,  \label{eq:hpnewconvderucp1}
\end{equation}
(\ref{Estrong}) holds with
\begin{equation*}
\mathcal{B}^{S}\left( u\right) _{s}=\int_{0}^{s}h\left( r,S_{r}\right)
dr+\int_{0}^{s}\left\langle \partial _{x}u\left( r,S_{r}\right) b_{1}\left(
r,S_{r}\right) \right\rangle dr-\int_{0}^{s}\left\langle \partial
_{x}u\left( r,S_{r}\right) b\left( r,S_{r}\right) \right\rangle dr.
\end{equation*}
\end{Corollary}

\textbf{Proof.} The result follows setting
\begin{equation*}
A_s = \int_0^s b_1(r,S_r) dr
\end{equation*}
in previous Theorem \ref{th:Itostrong}. \hfill
\hbox{\hskip 6pt\vrule width6pt height7pt
depth1pt  \hskip1pt}\bigskip

The above result depends on the extra assumption (\ref{eq:hpnewconvderucp1})
which is essential but not easy to check. We give first a special (but
useful) case where it holds and then an improvement for the nondegenerate
case. We have the following.

\begin{Remark}
\label{rm:suffperconvder}If
\begin{equation*}
\lim_{n\rightarrow +\infty }\partial _{x}u_{n}=\partial _{x}u\qquad \text{in
}C^{0}\left( \left[ 0,T\right] \times \mathbb{R}^{n}\right)
\end{equation*}
then Assumption (\ref{eq:hpnewconvderucp}) is verified. This means that the
result of Proposition \ref{pr:FRnew} above applies if we know that $u$ is a
strong solution in a more restrictive sense, i.e. substituting the point 2
of Definition \ref{df:solstrong} with
\begin{equation*}
\begin{array}{l}
u_{n}\longrightarrow u\text{ in }C^{0}\left( \left[ 0,T\right] \times
\mathbb{R}^{n}\right) \\
\partial _{x}u_{n}\longrightarrow \partial _{x}u\text{ in }C^{0}\left( \left[
0,T\right] \times \mathbb{R}^{n}\right) \\
h_{n}\longrightarrow h\text{ in }C^{0}\left( \left[ 0,T\right] \times
\mathbb{R}^{n}\right) \\
\phi _{n}\longrightarrow \phi \text{ in }C^{0}\left( \mathbb{R}^{n}\right) .
\end{array}
\end{equation*}
This is a particular case of our setting and it is the one used e.g
in \cite {CPDE,JMAA,EFJDE} to get the verification result. We can
say that in these works a result like Theorem \ref{th:Itostrong} is
proved under the assumption that $u$ is a strong solution in this
more restrictive sense. It is worth to note that in such simplified
setting, the proof of Theorem \ref {th:Itostrong} follows simply by
using standard convergence arguments. In particular there one does
not need to use the Fukushima-Dirichlet decomposition presented in
Section \ref{ITO}. So, from the methodological point of view there
is a serious difference with the result of Theorem \ref
{th:Itostrong}, see \cite{GR2}, Section 8
for comments.\hfill
\hbox{\hskip 6pt\vrule width6pt height7pt
depth1pt  \hskip1pt}\bigskip
\end{Remark}

A more significant achievement concerns the nondegenerate case. It is
illustrated in the following corollary.

\begin{Corollary}
\label{CGirsanov} \label{cr:Girsanov}We make the same assumption of
Corollary \ref{cor:ITOstrong} except (\ref{eq:hpnewconvderucp}) which we
replace by the assumption that
\begin{equation}
\sigma^{-1}\left( b_{1}-b\right) \text{ is bounded.}  \label{eq:hpGirsanov}
\end{equation}
where $\sigma^{-1}$ stands for the pseudo-inverse of $\sigma $. Then the
same conclusion of Corollary \ref{cor:ITOstrong} holds.
\end{Corollary}

\textbf{Proof. }Setting $t=0$ for simplicity we write
\begin{equation*}
\beta _{s}=W_{s}+\int_{0}^{s}\left( \sigma ^{-1}\left( b_{1}-b\right)
\right) \left( r,S_{r}\right) dr
\end{equation*}
and, applying Girsanov Theorem, there is a probability $\mathbb{Q}$
equivalent to $\mathbb{P}$ on $\left( \Omega ,\mathcal{F}_{T}\right)
$ such that $\left( \beta _{s}\right) _{s\in \left[ 0,T\right] }$ is
an $\left( \mathcal{F}_{s}\right) $-standard Brownian motion.

So, under $\mathbb{Q}$ the process $\left( S_{s}\right) $ fulfills the
equation
\begin{equation*}
dS_{s}=b\left( s,S_{s}\right) ds+\sigma \left( s,S_{s}\right) dW_{s};\qquad
S_{0}=x.
\end{equation*}
Under $\mathbb{Q}$, Assumption (\ref{eq:hpnewconvderucp}) is trivially
verified since $b=b_{1}$ and so we have (the tilde stands for operators
under the new probability $\mathbb{Q}$)
\begin{equation*}
\mathcal{\tilde{B}}^{S}\left( u\right) _{s}=\int_{0}^{s}h\left(
r,S_{r}\right) dr
\end{equation*}
and
\begin{equation*}
u\left( s,S_{s}\right) =u\left( 0,S_{0}\right) +\int_{0}^{s}\partial
_{x}u\left( r,S_{r}\right) \sigma \left( r,S_{r}\right) d\beta _{r}+\mathcal{%
\tilde{B}}^{S}\left( u\right) _{s}.
\end{equation*}
Expressing $\beta $ in terms of $W$ we obtain the result.\hfill
\hbox{\hskip 6pt\vrule width6pt height7pt
depth1pt  \hskip1pt}\bigskip

\subsection{Some useful consequences}

The previous results have some important consequences.

\begin{Remark}
\label{rm:FDdec2bis}From Remark \ref{rm:FDdec2} it follows that the
conclusion of the above 
Corollary \ref{cr:Girsanov} holds also if we stop the processes $\mathcal{B}%
(u)$ with a stopping time $t\leq \tau \leq T$. More precisely we have
\begin{equation*}
\mathcal{B}^{S}\left( u\right) _{s\wedge \tau }=\int_{0}^{s\wedge \tau
}h\left( r,S_{r}\right) dr-\int_{0}^{s\wedge \tau }\left\langle \partial
_{x}u\left( r,S_{r}\right) ,b_{1}\left( r,S_{r}\right) -b\left(
r,S_{r}\right) \right\rangle dr.
\end{equation*}
This fact will be useful in \cite{GR2}, Section 6. 
\hfill \hbox{\hskip 6pt\vrule width6pt height7pt depth1pt
\hskip1pt}\bigskip
\end{Remark}

\begin{Remark}
\label{rm:singat0representation}The results of the Corollary \ref
{cor:ITOstrong} and of Corollary \ref{cr:Girsanov} above still hold true
with suitable modifications if we assume that, instead of having $u\in
C^{0,1}\left( \left[ 0,T\right] \times \mathbb{R}^{n}\right) $:

\begin{itemize}
\item[(i)]  the strong solution $u$ belongs to $C^{0}\left( \left[ 0,T\right]
\times \mathbb{R}^{n}\right) \cap C^{0,1}\left( \left[ \varepsilon ,T\right]
\times \mathbb{R}^{n}\right) $ for every small $\varepsilon >0$;

\item[(ii)]  for some $\beta \in \left( 0,1\right) $ the map $\left(
t,x\right) \rightarrow t^{\beta }\partial _{x}u\left( t,x\right) $ belong to
$C^{0,1}\left( \left[ 0,T\right] \times \mathbb{R}^{n}\right) $.
\end{itemize}

The proof of Corollary \ref{cor:ITOstrong} in this case is a
straightforward generalization of the one presented above: we do not
give it here to avoid technicalities. In fact in proving the
verification theorem in \cite{GR2}, Section 6, we will deal with
initial data that are only continuous so with solutions $u$
satisfying (i) and (ii) above and possibly not $C^{0,1}$. This
difficulty will be faced directly in the proof of verification
Theorem 6.19 in \cite {GR2} by approximating the initial data with
$C^{1}$ ones, using (\ref {eq:reprmain}) and passing to the limit.
\hfill \hbox{\hskip 6pt\vrule width6pt height7pt depth1pt
\hskip1pt}\bigskip
\end{Remark}

\begin{Remark}
\label{rm:esempioestensione} From the above Corollary
\ref{cor:ITOstrong} it follows that the process
$\mathcal{B}^{S}\left( u\right) _{s}$ is in fact a semimartingale
(and also absolutely continuous). \hfill \hbox{\hskip 6pt\vrule
width6pt height7pt depth1pt \hskip1pt}\bigskip
\end{Remark}

\subsection{The elliptic case}

We devote the last part of this subsection to apply the same setting above
to elliptic problems. Consider the inhomogenous elliptic problem
\begin{equation}
\lambda u\left( x\right) +L_{0}u\left( x\right) +h\left( x\right) =0,\quad
\quad \forall x\in \mathbb{R}^{n}.  \label{eq:elllinear}
\end{equation}
where $D\left( L_{0}\right) =C^{2}\left( \mathbb{R}^{n}\right) $ and
\begin{equation*}
L_{0}u\left( x\right) =\left\langle b\left( x\right) ,\partial _{x}u\left(
s,x\right) \right\rangle +\frac{1}{2}\text{Tr }\left[ \sigma ^{\ast }\left(
x\right) \partial _{xx}u\left( x\right) \sigma \left( x\right) \right] .
\end{equation*}

\begin{Definition}
We say that $u$ is a strict solution to the elliptic problem (\ref
{eq:elllinear}) if $u\in D\left( L_{0}\right) $ and (\ref{eq:elllinear})
holds.
\end{Definition}

\begin{Definition}
\label{df:solstrongell}We say that $u$ is a strong solution to the elliptic
problem (\ref{eq:elllinear}) if there exists a sequence $\left( u_{n}\right)
\subseteq D\left( L_{0}\right) $ and a sequence $\left( h_{n}\right)
\subseteq C^{0}\left( \mathbb{R}^{n}\right) $, such that

\begin{enumerate}
\item  For every $n\in \mathbb{N}$ $u_{n}$ is a strict solution of the
problem
\begin{equation*}
\lambda u_{n}\left( x\right) -L_{0}u_{n}\left( x\right) =h_{n}\left(
x\right) ,\qquad \quad \quad \forall x\in \mathbb{R}^{n}.
\end{equation*}

\item  The following limits hold
\begin{equation*}
\begin{array}{l}
u_{n}\longrightarrow u\text{ in }C^{0}\left( \mathbb{R}^{n}\right) , \\
h_{n}\longrightarrow h\text{ in }C^{0}\left( \mathbb{R}^{n}\right) .
\end{array}
\end{equation*}
\end{enumerate}
\end{Definition}

Note that if $L$ is the closure of $L_{0}$ in $C^{0}\left( \mathbb{R}%
^{n}\right) $ then a strong solution $u$, by construction, belongs to $%
D\left( L\right) $. 
We now exploit the above setting to show that, for functions $u\in
C^{1}\left( \mathbb{R}^{n}\right) $ that are strong solutions of the elliptic
problem (\ref{eq:elllinear}) the following results holds; we omit the proof
as it is completely similar (and even simpler) to the one of Theorem \ref
{th:Itostrong} for the parabolic case.

\begin{Theorem}
\label{th:ITOstrongell} Let
\begin{equation*}
b_{1}:\Omega \times \mathbb{R}^{n}\rightarrow \mathbb{R}^{n},
\end{equation*}
be a continuous progressively measurable process (continuous in $x$) and
\begin{equation*}
b:\mathbb{R}^{n}\rightarrow \mathbb{R}^{n},\qquad \sigma :\mathbb{R}%
^{n}\rightarrow \mathcal{L}\left( \mathbb{R}^{m},\mathbb{R}^{n}\right)
\end{equation*}
be continuous functions. Let $S_{s}$ be a solution to the SDE
\begin{equation*}
dS_{s}=b_{1}\left( S_{s}\right) ds+\sigma \left( S_{s}\right) dW_{s};\qquad
S_{0}=x.
\end{equation*}
Let $u\in C^{1}\left( \mathbb{R}^{n}\right) $ be a strong solution of the
elliptic problem (\ref{eq:elllinear}). Assume that
\begin{equation}
\lim_{n\rightarrow +\infty }\int_{0}^{s}\left\langle \partial
_{x}u_{n}\left( S_{r}\right) -\partial _{x}u\left( S_{r}\right) ,b_{1}\left(
S_{r}\right) -b\left( S_{r}\right) \right\rangle dr=0,\qquad u.c.p.,
\label{eq:hpnewell}
\end{equation}
or that (\ref{eq:hpGirsanov}) holds. Then we have
\begin{equation}
u(S_{s})=u(S_{t})+\int_{t}^{s}\partial _{x}u\left( S_{r}\right) \sigma
(S_{r})dW_{r}+\mathcal{B}^{S}\left( u\right) _{s}-\mathcal{B}^{S}\left(
u\right) _{t},  \label{Estrongell}
\end{equation}
where
\begin{equation*}
\mathcal{B}^{S}\left( u\right) _{s}=\int_{0}^{s}h\left( S_{r}\right)
dr+\int_{0}^{s}\left\langle \partial _{x}u\left( S_{r}\right) ,b_{1}\left(
S_{r}\right) -b\left( S_{r}\right) \right\rangle dr.
\end{equation*}
\end{Theorem}

\bigskip

\textbf{ACKNOWLEDGEMENTS} The authors wish to thank the Referee and the
Editor for the careful reading and for the motivating suggestions.

\end{document}